\documentclass{amsart}

\usepackage{fancyhdr}
\usepackage{lastpage}
\usepackage{yhmath}

\newcommand{\bdis}{\begin{displaymath}}
\newcommand{\edis}{\end{displaymath}}
\newcommand{\be}{\begin{equation}}
\newcommand{\ee}{\end{equation}}

\newcommand{\mbb}{\mathbb}
\newcommand{\mcal}{\mathcal}

\newcommand{\vp}{\varphi}
\newcommand{\vt}{\vartheta}

\newcommand{\zf}{\zeta\left(\frac{1}{2}+it\right)}

\DeclareMathOperator{\re}{Re}

\theoremstyle{definition}

\theoremstyle{remark}

\newtheorem*{mydef1}{{\bf Theorem}}

\newtheorem*{mydef2}{{\bf Definition}}

\newtheorem*{mydef5A}{{\bf Lemma $\alpha$}}
\newtheorem*{mydef5B}{{\bf Lemma $\beta$}}

\newtheorem*{mydef8}{{\bf Remark}}

\numberwithin{equation}{section}

\begin{document}

\title{Improvement of the theorem of Hardy-Littlewood on density of zeros of the function $\zf$}

\author{Jan Moser}

\address{Department of Mathematical Analysis and Numerical Mathematics, Comenius University, Mlynska Dolina M105, 842 48 Bratislava, SLOVAKIA}

\email{jan.mozer@fmph.uniba.sk}


\begin{abstract}
In this paper we improve classical Hardy-Littlewood exponent $1/2$ by about  $16.6\%$ 62 years after the original
result. This result is the first step to prove the Selberg's hypothesis (1942). In order to reach our purpose we
use discrete method. This paper is the English version of our paper \cite{4}.
\end{abstract}

\maketitle

\section{Introduction}

Let us remind that Hardy and Littlewood have proved (see \cite{1}, p. 283) the following classical
estimate

\be \label{1.1}
N_0(T+T^{1/2+\epsilon}-N_0(T)>A(\epsilon)T^{1/2+\epsilon},\quad T\geq T_0(\epsilon)
\ee
for every fixed small positive $\epsilon$, where $N_0(T)$ denotes the number of zeros of the function
\bdis
\zf,\quad t\in (0,T],
\edis
and $A(\epsilon)$ is a constant depending on the choice of $\epsilon$.

In this paper we shall prove, by a discrete method, a theorem on number of good segments (definition of these
is placed below). Next, as a Corollary, we obtain the following estimate
\be \label{1.2}
N_0(T+T^{5/12}\psi\ln^3T)-N_0(T)>A(\psi)T^{5/12}\psi\ln^3T,
\ee
where $\psi=\psi(T)$ means an arbitrarily slowly increasing function unbounded from above, for example
\bdis
\psi=\ln\ln\dots\ln T,
\edis
and $A(\psi)$ is a constant depending upon the choice of $\psi$.

Next, let us remind that A. Selberg (see \cite{5}, p. 46, Theorem A) has obtained the fundamental result by means
of the estimate
\bdis
N_0(T+T^{1/2+\epsilon})-N_0(T)>A(\epsilon)T^{1/2+\epsilon}\ln T
\edis
and he has raised the following hypothesis
\be \label{1.3}
\frac 12 \longrightarrow a:\ a<\frac 12.
\ee

\begin{mydef8}
We notice explicitly that
\begin{itemize}
 \item[(a)] our improvement (\ref{1.2})
 \bdis
 \frac 12 \longrightarrow \frac{5}{12}
 \edis
 of the Hardy-Littlewood exponent in (\ref{1.1}) represents the part $16.6\%$,
 \item[(b)] estimate (\ref{1.2}) is the first step on a way to prove the Selberg's hypothesis (\ref{1.3}).
\end{itemize}
\end{mydef8}

\section{Theorem}

Let (see \cite{7}, pp. 79, 329)
\be \label{2.1}
Z(t)=e^{i\vt(t)}\zf ,
\ee
\be \label{2.2}
\vt(t)=-\frac t2\ln\pi+\text{Im}\ln\Gamma\left(\frac 12+i\frac t2\right)=\vt_1(t)+\mcal{O}\left(\frac 1t\right),
\ee
\be \label{2.3}
\vp_1(t)=\frac t2\ln\frac{t}{2\pi}-\frac t2-\frac{\pi}{8},
\ee
and $\{\bar{t}_\nu\}$ denote the sequence defined by the following formula (see \cite{3}, (2))
\be \label{2.4}
\vt_1(\bar{t}_\nu)=\frac{\pi}{2}\nu,\ \nu=1,2,\dots
\ee
Next, let
\be \label{2.5}
\omega=\frac{\pi}{\ln\frac{T}{2\pi}},\ U=T^{5/12}\psi\ln^3T,\ M_1=M_1(\delta,T)=[\delta\ln T],\ \delta>1.
\ee

\begin{mydef2}
We shall call the segment
\be \label{2.6}
[\bar{t}_\nu+k(\nu)\omega,\bar{t}_\nu+(k(\nu)+1)\omega],
\ee
where
\be \label{2.7}
\bar{t}_\nu\in [T,T+U],\ 0\leq k(\nu)\leq M_1 ,
\ee
and $k(\nu)\in \mbb{N}_0$, as \emph{the good segment} (comp. \cite{6}, \cite{2}) if
\be \label{2.8}
Z(\bar{t}_\nu+k(\nu)\omega)\cdot Z(\bar{t}_\nu+(k(\nu)+1)\omega)<0.
\ee
\end{mydef2}

Let $G(T,U,\delta)$ denote the number of non-intersecting good segments within the segment $[T,T+U]$. The following
theorem holds true.

\begin{mydef1}
There are
\bdis
\delta_0>0,\ A(\psi,\delta_0)>0,\ T_0(\psi,\delta_0)>0
\edis
such that
\be \label{2.9}
G(T,U,\delta_0)>A(\psi,\delta_0)U,\ T\geq T_0(\psi,\delta_0).
\ee
\end{mydef1}

Since the good segment contains a zero point of the odd order of the function
\bdis
\zf
\edis
(see (\ref{2.1}), (\ref{2.8})) then the estimate (\ref{1.2}) follows, where, of course, $A(\psi)=A(\psi,\delta_0)$.

\section{Main lemmas and proof of Theorem}

\subsection{}

Let (see \cite{3}, (3))
\be \label{3.1}
J=\sum_{k=0}^M \sum_{l=0}^M\sum_{T\leq \bar{t}_\nu\leq T+U} Z(\bar{t}_\nu+k\omega)Z(\bar{t}_\nu+l\omega).
\ee
We have (see Theorem from \cite{3}) the following

\begin{mydef5A}
\be \label{3.2}
J=AMU\ln^2T+o(MU\ln^2T),
\ee
where
\be \label{3.3}
U=T^{5/12}\psi\ln^3T,\ \ln T<M<\sqrt[3]{\psi}\ln T,
\ee
and $A>0$ is an absolute constant.
\end{mydef5A}

Next, we put
\be \label{3.4}
N=\sum_{T\leq \bar{t}_\nu\leq T+U} |K|^2,
\ee
where
\be \label{3.5}
K=\sum_{k=0}^M
\left\{ e^{-i\vt(\bar{t}_\nu+k\omega)}Z(\bar{t}_\nu+k\omega)-1\right\}.
\ee
The following lemma holds true

\begin{mydef5B}
\be \label{3.6}
N=\mcal{O}(MU\ln^2T).
\ee
\end{mydef5B}

\subsection{}

Now we use our main lemmas to complete the proof of the Theorem. Let
\be \label{3.7}
J(\bar{t}_\nu)=\sum_{k=0}^M Z(\bar{t}_\nu+k\omega),
\ee
\be \label{3.8}
L(\bar{t}_\nu)=\sum_{k=0}^M |Z(\bar{t}_\nu+k\omega)|,
\ee
and
\bdis
R=R(T,U)
\edis
denote the number of $\bar{t}_\nu^*$ with the property
\be \label{3.9}
\bar{t}_\nu^*\in [T,T+U] \ \Rightarrow \ |J(\bar{t}_\nu^*)|=L(\bar{t}_\nu^*).
\ee
It is clear that the members of the sequence
\be \label{3.10}
\{ Z(\bar{t}_\nu^*+k\omega)\}_{k=0}^M
\ee
preserve their signs. Thus
\be \label{3.11}
\sum_{\bar{t}_\nu^*}|J|=\sum_{\bar{t}_\nu^*}L.
\ee
Next, we have (see (\ref{3.2}), (\ref{3.5}))
\be \label{3.12}
\sum_{\bar{t}_\nu^*}|J|\leq \sqrt{R}\left(\sum_{\bar{t}_\nu^*}J^2\right)^{1/2}=\sqrt{RJ}<
A\sqrt{RMU\ln^2T},
\ee
\be \label{3.13}
\begin{split}
 & L=\sum_{k=0}^M |Z(\bar{t}_\nu+k\omega)|\geq \left|\sum_{k=0}^M e^{-i\vt(\bar{t}_\nu+k\omega)}
 Z(\bar{t}_\nu+k\omega)\right|=\\
 & = |K+M+1|\geq M+1-|K|.
\end{split}
\ee
Now we have (see (\ref{3.6}))
\be \label{3.14}
\begin{split}
 & \sum_{\bar{t}_\nu^*}L\geq (M+1)R-\sum_{\bar{t}_\nu^*}|K|\geq \\
 & \geq (M+1)R-\sqrt{R}\sqrt{\sum_{\bar{t}_\nu^*}|K|^2}\geq (M+1)R-\sqrt{RN}>\\
 & > (M+1)R-A\sqrt{RMU\ln^2T},
\end{split}
\ee
and (see (\ref{3.11}),(\ref{3.12}),(\ref{3.14}))
\be \label{3.15}
(M+1)R<A\sqrt{RMU\ln^2T}.
\ee
Consequently, we obtain
\be \label{3.16}
R<A\frac{U\ln^2T}{M}.
\ee
Now, we divide the number (comp. \cite{3}, (8))
\be \label{3.17}
Q_1=\frac 1\pi U\ln\frac{T}{2\pi}+\mcal{O}\left(\frac{U^2}{T}\right)
\ee
of values
\bdis
\bar{t}_\nu\in [T,T+U]
\edis
into
\bdis
\left[\frac{Q_1}{2M}\right]
\edis
pairs of abutting parts $j_1,j_2$, each except the last $j_2$, of length $M$ (comp. \cite{7}, p. 226). Let now
$\mu$ denote the number of parts $j_1$ consisting entirely of the points $\bar{t}_\nu^*$. Then by
(\ref{3.16}) we have
\be \label{3.18}
\mu M<A\frac{U\ln^2T}{M}\ \Rightarrow \ \mu<AU\left(\frac{\ln T}{M}\right)^2.
\ee
If (see (\ref{2.5}))
\bdis
M=M_1=[\delta\ln T],
\edis
then we obtain by (\ref{2.5}), (\ref{3.17}) and (\ref{3.18})
\be \label{3.19}
\begin{split}
 & \left[\frac{Q_1}{2M}\right]-\mu> A_1\frac{U\ln T}{M}-A_2\frac{U^2}{MT}-\mu>
 \frac{1}{\delta_0}\left( A_3-\frac{A_4}{\delta_0}\right)U-A_5> \\
 & > A(\psi,\delta_0)U
\end{split}
\ee
for sufficiently big $\delta$, ($=\delta_0$ say). This inequality gives the estimate from below
for thu number of parts $j_1$ such that every of these contains at least one point
$\bar{t}_\nu$ for which (see (\ref{3.7}) -- (\ref{3.9}))
\be \label{3.20}
|J(\bar{t}_\nu)|\not= L(\bar{t}_\nu).
\ee
Of course, for such a point we have that the members of the sequence
\bdis
\{ Z(\bar{t}_\nu+k\omega)\}_{k=0}^M
\edis
change sign, i.e. there is
\bdis
k(\nu)\in [0,M(\delta_0)]
\edis
such that (\ref{2.8}) holds true. Since there is at least (see (\ref{3.19})) $A(\psi,\delta_0)U$
of parts $j_1$ of this kind, then (\ref{2.9}) follows.

Proof of Lemma $\beta$ is in what follows.

\section{Decomposition of the sum $N$}

Let (comp. (\ref{3.5}))
\be \label{4.1}
\begin{split}
 & K_1=\re\left\{ e^{-i\vt(\bar{t}_\nu+k\omega)}Z(\bar{t}_\nu+k\omega)-1\right\}\cdot
 \left\{ e^{i\vt(\bar{t}_\nu+l\omega)}Z(\bar{t}_\nu+l\omega)-1\right\}= \\
 & = \re \left\{ (e^{-i\vt_k}Z_k-1)(e^{i\vt_l}Z_l-1)\right\}= \\
 & = Z_kZ_l\cos(\vt_k-\vt_l)-Z_k\cos\vt_k-Z_l\cos\vt_l+1.
\end{split}
\ee
Putting
\be \label{4.2}
Z_k=2\cos\vt_k+\bar{Z}_k,\ \bar{t}_k=\bar{t}_\nu+k\omega,
\ee
where (see \cite{3}, (117))
\be \label{4.3}
\bar{Z}_k=2\sum_{2\leq n<P_0}\frac{1}{\sqrt{n}}
\cos(\vt_k-\bar{t}_k\ln n)+\mcal{O}(T^{-1/4}),\ P_0=\sqrt{\frac{T}{2\pi}},
\ee
we obtain
\be \label{4.4}
\begin{split}
 & K_1=\bar{Z}_k\bar{Z}_l\cos(\vt_k-\vt_l)+\\
 & + 2\bar{Z}_k\cos\vt_l\cos(\vt_k-\vt_l)+2\bar{Z}_l\cos\vt_k\cos(\vt_k-\vt_l)-\\
 & - \bar{Z}_k\cos\vt_k-\bar{Z}_l\cos\vt_l+ \\
 & + 4\cos\vt_k\cos\vt_l\cos(\vt_k-\vt_l)-2\cos^2\vt_k-2\cos^2\vt_l+1.
\end{split}
\ee
Hence, (see (\ref{3.4}), (\ref{3.5}), (\ref{4.1}), (\ref{4.4}))
\be \label{4.5}
\begin{split}
 & N=\sum_{T\leq\bar{t}_\nu\leq T+U }\sum_{k=0}^M\sum_{l=0}^M K_1=
 \sum_{\bar{t}_\nu}\sum_k\sum_l \bar{Z}_k \bar{Z}_l \cos(\vt_k-\vt_l)+\\
 & + 4\sum_{\bar{t}_\nu}\sum_k\sum_l \bar{Z}_k\cos\vt_l\cos(\vt_k-\vt_l)-
 2\sum_{\bar{t}_\nu}\sum_k\sum_l\bar{Z}_k\cos\vt_k+ \\
 & + \sum_{\bar{t}_\nu}\sum_k\sum_l
 \{ 4\cos\vt_k\cos\vt_l\cos(\vt_k-\vt_l)-4\cos^2\vt_k+1\}= \\
 & = w_1+w_2+w_3+w_4.
\end{split}
\ee

\section{Estimate of $w_4$}

Next, we obtain, (see (\ref{2.2}) -- (\ref{2.4}) and \cite{3}, (118), (119)) that
\be \label{5.1}
\begin{split}
 & 4\cos\vt_k\cos\vt_l\cos(\vt_k-\vt_l)= \\
 & = 2\cos^2(\vt_k-\vt_l)+2\cos(\vt_k+\vt_l)\cos(\vt_k-\vt_l)= \\
 & = 1+\cos\{ 2(\vt_k-\vt_l)\}+\cos(2\vt_k)+\cos(2\vt_l)= \\
 & = 1+\cos\{2(\vt_{1,k}-\vt_{1,l})\}+\cos(2\vt_{1,k})+\cos(2\vt_{1,l})+
 \mcal{O}\left(\frac{1}{T}\right)= \\
 & = 1+\cos\{ 2(k-l)\omega\ln P_0\}+\cos(\pi\nu+2k\omega\ln P_0)+\\
 & + \cos(\pi\nu+2l\omega\ln P_0)+\mcal{O}\left(\frac{MU}{T\ln T}\right)+
 \mcal{O}\left(\frac{1}{T}\right)= \\
 & = 1+(-1)^{k+l}+(-1)^{\nu+k}+(-1)^{\nu+l}+\mcal{O}\left(\frac{MU}{T\ln T}\right),
\end{split}
\ee
\be \label{5.2}
\begin{split}
 & -4\cos^2\vt_k=-2-2\cos(2\vt_k)= \\
 & = -2-2\cos(\pi\nu+2k\omega\ln P_0)+\mcal{O}\left(\frac{MU}{T\ln T}\right)= \\
 & = -2-2(-1)^{\nu+k}+\mcal{O}\left(\frac{MU}{T\ln T}\right),
\end{split}
\ee
since (see (\ref{2.5}), (\ref{4.3}))
\be \label{5.3}
2\omega \ln P_0=\pi .
\ee
Finally, we have (see (\ref{3.17}), (\ref{4.5}), (\ref{5.1}), (\ref{5.2}))
\be \label{5.4}
\begin{split}
 & w_4=\sum_{\bar{t}_\nu}\sum_k\sum_l (-1)^{k+l}+\mcal{O}
 \left( M^2U\ln T\frac{MU}{T\ln T}\right)= \\
 & = \mcal{O}(U\ln T)+\mcal{O}\left(\frac{M^3U^2}{T}\right).
\end{split}
\ee

\section{Estimate of $w_1$}

Further, we have (see (\ref{2.2}, (\ref{2.3}), (\ref{4.3}) and \cite{3}, (119), (120))
\be \label{6.1}
\begin{split}
 & \bar{Z}_k\bar{Z}_l=2\sum_m\sum_n \frac{1}{\sqrt{mn}}
 \cos(\vt_k+\vt_l-\bar{t}_k\ln n-\bar{t}_l \ln m)+ \\
 & +2\sum_m\sum_n \frac{1}{\sqrt{mn}}
 \cos(\vt_k-\vt_l-\bar{t}_k\ln n+\bar{t}_l\ln m)+\mcal{O}(T^{-1/12}\ln T)= \\
 & = 2\sum_m\sum_n\frac{(-1)^\nu}{\sqrt{mn}}
 \cos\{ \bar{t}_\nu\ln(mn)-(k+l)\omega\ln P_0+k\omega\ln n+l\omega\ln m\}+ \\
 & + 2\sum_m\sum_n\frac{1}{\sqrt{mn}}
 \cos\left\{ \bar{t}_\nu\ln\frac nm+(l-k)\omega\ln P_0+k\omega\ln n-l\omega\ln m\right\}+ \\
 & + \mcal{O}\left(\frac{MU}{\sqrt{T}\ln T}\right)+\mcal{O}(T^{-1/12}\ln T).
\end{split}
\ee
Next, (see (\ref{6.1}), comp. \cite{3}, (121))
\be \label{6.2}
\begin{split}
 & \bar{Z}_k\bar{Z}_l\cos(\vt_k-\vt_l)=\bar{Z}_k\bar{Z}_l\cos\{(k-l)\omega\ln P_0\}+
 \mcal{O}\left(\frac{MU\ln T}{T^{2/3}}\right)= \\
 & = \sum_m\sum_n\frac{(-1)^\nu}{\sqrt{mn}}
 \cos\{ \bar{t}_\nu\ln(mn)-2l\omega\ln P_0+k\omega\ln n+l\omega\ln m\}+ \\
 & + \sum_m\sum_n\frac{(-1)^\nu}{\sqrt{mn}}
  \cos\{ \bar{t}_\nu\ln(mn)-2k\omega\ln P_0+k\omega\ln n+l\omega\ln m\}+ \\
 & + \sum_m\sum_n\frac{1}{\sqrt{mn}}
 \cos\left( \bar{t}_\nu\ln\frac nm+k\omega\ln n-l\omega\ln m\right)+ \\
 & + \sum_m\sum_n\frac{1}{\sqrt{mn}}
 \cos\left\{ \bar{t}_\nu\ln\frac nm-2(k-l)\omega\ln P_0+k\omega\ln n-l\omega\ln m\right\}+ \\
 & + \mcal{O}\left(\frac{MU}{\sqrt{T}\ln T}\right)+\mcal{O}(T^{-1/2}\ln T)+
 \mcal{O}\left(\frac{MU\ln T}{T^{2/3}}\right)= \\
 & = s_5+s_6+s_7+s_8+\mcal{O}\left(\frac{MU}{\sqrt{T}\ln T}\right)+
 \mcal{O}(T^{-1/12}\ln T).
\end{split}
\ee
First of all
\be \label{6.3}
\sum_k\sum_l\sum_{\bar{t}_\nu}(s_5+s_6+s_7^{(m\not= n)}+s_8^{(m\not=n)})=
\mcal{O}(M^3T^{5/12}\ln^3T)
\ee
by lemmas of type B and C from \cite{3}. Next, (comp. \cite{3}, (94), (104))
\be \label{6.4}
s_{71}=\sum_k\sum_l s_7^{(m\not= n)}=\sum_{2\leq n<P_0}\frac 1n G(M+1,\omega\ln n).
\ee
Since (see (\ref{2.5}))
\bdis
0<\frac 12 \omega\ln n=\frac{\pi}{4}\frac{\ln n}{\ln P_0}<\frac{\pi}{4},
\edis
then
\be \label{6.5}
\sin\left(\frac 12\omega\ln n\right)>A\omega\ln n,
\ee
and
\be \label{6.6}
s_{71}=\mcal{O}\left(\frac{1}{\omega^2}\sum_{n=2}^\infty \frac{1}{n\ln^2n}\right)=
\mcal{O}(\ln^3T).
\ee
Consequently, (see (\ref{3.17}))
\be \label{6.7}
\sum_{\bar{t}_\nu} s_{71}=\mcal{O}(U\ln T\ln^2T)=\mcal{O}(U\ln^3T).
\ee
Next we obtain by similar way that
\be \label{6.8}
s_{81}=\sum_k\sum_l s_8^{(m=n)}=\sum_{2\leq n<P_0}\frac 1n
G\left( M+1,\omega\ln\frac{P_0^2}{n}\right).
\ee
Since
\bdis
\frac 12\omega\ln\frac{P_0^2}{n}=\frac{\pi}{2\ln P_0^2}\ln\frac{P_0^2}{n},
\edis
then
\bdis
\frac{\pi}{4}<\frac 12\omega\ln\frac{P_0^2}{n}<\frac{\pi}{2},
\edis
and
\be \label{6.9}
\sin\left(\frac 12\omega\ln\frac{P_0^2}{n}\right)>A>0.
\ee
Hence
\be \label{6.10}
s_{81}=\mcal{O}\left(\sum_{2\leq n<P_0}\frac 1n\right)=\mcal{O}(\ln T),
\ee
and consequently,
\be \label{6.11}
\sum_{\bar{t}_\nu} s_{81}=\mcal{O}(U\ln^2T).
\ee
Finally, we obtain (see (\ref{4.5}), (\ref{6.2}), (\ref{6.3}), (\ref{6.7}), (\ref{6.11})) in the case
(\ref{3.3}) that
\be \label{6.12}
\begin{split}
 & w_1=\mcal{O}(M^3T^{5/12}\ln^3T)+\mcal{O}(U\ln^3T)+\\
 & + \mcal{O}\left( \frac{M^3U^2}{\sqrt{T}}\right)+\mcal{O}(M^2UT^{-1/12}\ln^3T)=
 \mcal{O}(MU\ln^2T).
\end{split}
\ee

\section{Estimates of $w_2,w_3$}

First of all, we have (see (\ref{4.2}), (\ref{4.3}), comp. (\ref{4.5}))
\be \label{7.1}
\begin{split}
 & 4\bar{Z}_k\cos\vt_k\cos(\vt_k-\vt_l)= \\
 & = \sum_n \frac{1}{\sqrt{n}}\cos(\bar{t}_\nu\ln n)+\sum_{n}\frac{1}{\sqrt{n}}\cos(2\vt_k-\bar{t}_\nu\ln n)+ \\
 & + \sum_n \frac{1}{\sqrt{n}}\cos(2\vt_l-\bar{t}_k\ln n)+
 \sum_n \frac{1}{\sqrt{n}}\cos(2\vt_k-2\vt_l-\bar{t}_k\ln n)+ \\
 & + \mcal{O}(T^{-1/4})= \\
 & = w_{21}+w_{22}+w_{23}+w_{24}+\mcal{O}(T^{-1/4}),\quad \bar{t}_k=\bar{t}_\nu+k\omega.
\end{split}
\ee
Next, we have (see (\ref{2.2}), (\ref{2.3}), (\ref{3.17}), (\ref{5.3}) and \cite{3}, (118))
\be \label{7.2}
\begin{split}
 & w_{221}=\sum_{T\leq\bar{t}_\nu\leq T+U}w_{22}= \\
 & = \sum_n \frac{1}{\sqrt{n}}\sum_{\bar{t}_\nu}\cos(\pi\nu+2k\omega\ln P_0-\bar{t}_k\ln n)+ \\
 & +
 \mcal{O}\left\{ U\ln T\cdot T^{1/4}\left(\frac 1T+\frac{MU}{T\ln T}\right)\right\}= \\
 & = (-1)^k\sum_n \frac{1}{\sqrt{n}}\sum_{\bar{t}_\nu}
 \cos(\pi\nu-\bar{t}_\nu\ln n-k\omega\ln n)+\mcal{O}(MU^2T^{-1/4})= \\
 & = (-1)^k\sum_n \frac{1}{\sqrt{n}}\sum_{\bar{t}_\nu}
 \cos(\pi\nu+\bar{t}_\nu\ln n+k\omega\ln n)+\mcal{O}(MU^2T^{-1/4}).
\end{split}
\ee
The following function corresponds to the $\bar{t}_\nu$-sum (comp. \cite{6}, pp. 99, 100)
\bdis
\chi(\nu)=\frac{1}{2\pi}(\pi\nu+\bar{t}_\nu+k\omega\ln n),\ \bar{t}_\nu\in [T,T+U].
\edis
Consequently, (see (\ref{2.3}), (\ref{2.4}))
\bdis
\begin{split}
 & \chi'(\nu)=\frac 12+\frac{1}{2\pi}\frac{\pi}{2}\frac{\ln n}{\ln\sqrt{\frac{\bar{t}_\nu}{2\pi}}}=
 \frac 12+\frac 14\frac{\ln n}{\ln P_0}+\mcal{O}\left(\frac{U}{T\ln T}\right), \\
 & \chi''(\nu)<0.
\end{split}
\edis
Hence
\bdis
\chi'(\nu)\in \left(\frac 12,\frac 34+\epsilon\right),\ \epsilon\in (0,1/4),
\edis
and (see (\ref{7.2}), comp. \cite{6}, p. 100)
\bdis
w_{221}=\mcal{O}\left(\sum_n \frac{1}{\sqrt{n}}\right)+\mcal{O}(MU^2T^{-3/4})=\mcal{O}(T^{1/4}).
\edis
Similar estimates can be obtained also for
\bdis
w_{211}, w_{231}, w_{241}.
\edis
Consequently, (see (\ref{4.5}), (\ref{7.1}))
\be \label{7.3}
w_2=\mcal{O}(M^2T^{1/4}),
\ee
and similarly we obtain the estimate
\be \label{7.4}
w_3=\mcal{O}(M^2T^{1/4}).
\ee
Finally, we obtain by (\ref{2.5}), (\ref{4.5}), (\ref{5.4}), (\ref{6.12}), (\ref{7.3}) and
(\ref{7.4}) that
\bdis
N=\mcal{O}(MU\ln^2T),
\edis
i. e. the estimate (\ref{3.6}) holds true.

\thanks{I would like to thank Michal Demetrian for helping me with the electronic version of this work.}

\end{document}